\documentclass{ifacconf}

\usepackage{graphicx}      
\usepackage{natbib}        
\usepackage{amsmath,amssymb}

\newcommand{\R}{{\mathbb R}}

\newcommand{\bk}[2]{\left\langle #1,#2\right\rangle}

\newcommand{\cL}{{\mathcal L}}
\newcommand{\cH}{{\mathcal H}}

\newcommand{\cA}{{\mathcal A}}

\newcommand{\bE}{{\mathbf E}}

\newcommand{\one}{{\mathbf 1}}

\newcommand{\eps}{\epsilon}

\newcommand{\wrt}{with respect to }

\begin{document}
	\begin{frontmatter}
		
		\title{Reachability Analysis of Randomly Perturbed Hamiltonian Systems} 

		\author{Carsten Hartmann, Lara Neureither, Markus Strehlau}

		\address{Institut f\"ur Mathematik, Brandenburgische Technische Universit\"at Cottbus-Senftenberg, Germany\\ (e-mail: \{carsten.hartmann,lara.neureither,markus.strehlau\}@b-tu.de)}

		\begin{abstract}                
			In this paper, we revisit energy-based concepts of controllability and reformulate them for control-affine nonlinear systems perturbed by white noise. Specifically, we discuss the relation between controllability of deterministic systems and the corresponding stochastic control systems in the limit of small noise and in the case in which the target state is a measurable subset of the state space.  We derive computable expressions for hitting probabilities and mean first hitting times in terms of empirical Gramians, when the dynamics is given by a Hamiltonian system perturbed by dissipation and noise, and  provide an easily computable expression for the corresponding controllability function as a function of a subset of the state variables. 
		\end{abstract}
		
		
	\end{frontmatter}

	\section{Introduction}
	
	Since the seminal work of \cite{moore1981}, energy-based controllability and reachability analysis has been playing a key role in model order reduction of linear control systems. The extension to nonlinear systems has been -- and still is -- an active field of research that comprises both empirical and algebraic approaches; see \cite{Condon2005,Hahn2002,lall2002,Scherpen1993,scherpen2021}. 
	
	While it is possible to put the controllability of linear and nonlinear systems under the umbrella of a geometric theory, going back to the influential work by \cite{krener1974}, the energy-based approach is not free of ambiguities as has been pointed out by several authors, e.g. \cite{benner2011,gray1998}. 
	
	Here we take side with energy-based concepts, but we adopt an alternative perspective, using that complete controllability and reachability of control-affine systems are intimately related to ergodicity of stochastic differential equations (see e.g. \cite{mattingly2002}). This connection allows for computing either empirical controllability functions or empirical controllability Gramians associated with a control system by Monte Carlo; see  \cite{Newman1999,hartmann2013,kashima2016}. Even though our approach relies on white noise inputs that can be interpreted as signals that are uniform in the frequency domain, the interpretation as a control is subtle, since such a control would have infinite $L^2$-norm. (The same goes for Dirac-valued inputs that are typically used to compute empirical Gramians based on impulse responses.) So one might ask why considering white noise inputs is a sensible thing to do, and it turns out that the relevance of white noise as forcing of a control system lies in its regularising effect on the controllability function, in that it mimics any piecewise linear control in the small noise regime where the approximation is in the mean square sense (see e.g. \cite{supportThm}). Hence, from the point of view of reachability analysis, we can think of a small-noise diffusion as a (deterministic) control system with square integrable controls. White noise inputs are also employed in the analysis of obstacle avoidance problems, when the quantity of interest is the probability to reach a set before some finite time $T$ or the probability to reach one set before another (see e.g. \cite{abate2008,Bujorianu2012}).

	In this paper we study (a) the connection between the energy-based reachability concept for control-affine nonlinear systems of \cite{Scherpen1993} and the stochastic set reachability (i.e. hitting probability) problem, (b) discuss some specifics of mechanical or port-Hamiltonian-like systems in the small noise regime (following our own work  \cite{HartmannSchuette2008,neureither2017time}), and (c) discuss the connection of the stochastic reachability problem with coarse-graining and its implication for the computation of empirical Gramians and controllability functions.   
	
	The paper is outlined as follows: Set reachability and finite-time controllability of deterministic and stochastic systems are discussed in Sections \ref{sec:reachability} and \ref{sec:SOC}. Section \ref{sec:smallnoise} is devoted to the question of how to compute the controllability function of small noise diffusions, with a particular focus on linear and nonlinear port-Hamiltonian-like systems. The theoretical results are illustrated  in Section \ref{sec:num} with a numerical study of a stochastically perturbed double pendulum. We summarize our findings in Section \ref{sec:sum}.

	\section{Controllability function of a nonlinear control system}\label{sec:reachability}
	
	To begin with, we consider an affine control system 
	\begin{equation}\label{ode}
		\frac{d}{ds}x(s) = f(x(s)) + g(x(s)) u(s)\,,\quad x(t)=x\,
	\end{equation}
	in $\R^n$ for a continuous control $u\colon\R\to M$ taking values in a set $M\subset \R^m$. 
	Here $(x,t)\in(\R^n,\R)$ denotes the initial data, and we assume that $f\colon\R^n\to\R^n$ and $g\colon\R^n\to\R^{n\times m}$ are smooth and satisfy suitable growth conditions. 
	
	A typical control task consists in finding a, say, continuous and bounded or square-integrable control $u$ that allows to reach a given target state $y\in\R^n$. Under certain assumptions, the set of all reachable states $y\in\R^n$ can be characterised by the quadratic cost functional 
	\begin{equation}\label{scherpenCost}
		I(u) = \frac{1}{2}\int_{-\infty}^{0} |u(s)|^2 ds
	\end{equation}
	having a finite value under the constraint
	\begin{equation}\label{scherpenConstraint}
		\lim_{t\to-\infty}x(t)=0\,,\quad x(0)=y \,.
	\end{equation}
	The following result is due to \cite{Scherpen1993} and establishes the relation o the controllability function with the value function of an optimal control problem; see also \cite{Newman1999}:
	\begin{thm}[\cite{Scherpen1993}]
		Let $L\colon\R^n\to[0,\infty)$ be the value function associated with (\ref{scherpenCost})--(\ref{scherpenConstraint}), considered as a function of the target state $y$. 
		Then there exists a neighbourhood $N\subset\R^n$ of $0$, such that $L$ is the unique viscosity solution of the dynamic programming equation
		\begin{equation}\label{scherpenHJB}
			f\cdot\nabla L + \frac{1}{2}|\nabla L|^2_{gg^T}=0\,,\quad L(0)=0\,,
		\end{equation}
		in $N\subset\R^n$, such that $0\in\R^n$ with $f(0)=0$ is an asymptotically stable equilibrium of 
		\begin{equation}\label{scherpenStable}
			f_{-}^* = -(f+gg^T\nabla L)\,.
		\end{equation} 
		Here $|\cdot|_{gg^T}$ denotes the weighted Euclidean semi-norm defined by $|z|_{Q}^2=z^TQz$ for a symmetric and positive semi-definite matrix $Q=Q^T\ge 0$. 
	\end{thm}
	
	\subsection{Computing the controllability function}
	
	The reachable states, for which the controllability function is finite, are often the basis of model reduction schemes for control systems. However there are only few cases in which an explicit solution to (\ref{scherpenHJB})--(\ref{scherpenStable}) is known, beyond the linear case in which solving the dynamic programming PDE can be reduced to solving an algebraic matrix equation of Lyapunov type. One such example of a nonlinear system, in which an explicit solution of (\ref{scherpenHJB})--(\ref{scherpenStable}) is known is 
	\begin{equation}\label{phs}
		\frac{d}{dt}x(s) = (J-D)\nabla H(x(s)) + \sqrt{2 D} u(s)\,,
	\end{equation}
	in which case 
	\begin{equation}\label{V-H}
		L(x) = H(x) - H(0).
	\end{equation}
	Here $H\colon\R^{n}\to\R$ is a smooth Hamiltonian that is bounded from below and at least quadratically growing at infinity, and $J=-J^T$ and $D=D^T\ge 0$ are constant $n\times n$ matrices where $n=2d$ is even. 
	%
	It can be readily verified that (\ref{V-H}) is a solution to (\ref{scherpenHJB})--(\ref{scherpenStable}), with $f = (J-D)\nabla H$ and $gg^T = 2 D$. 
	
	An alternative to the explicit solution is the approximation of the controllability function by Monte Carlo. The latter  is related to the invariant measure $\mu$ associated with the It\^{o} stochastic differential equation (SDE)
	\begin{equation}\label{sde}
		dX_s = f(X_s)ds + \sqrt{\eps} g(X_s)dW_s\,,
	\end{equation}
	for some $\eps>0$, where $W=(W_s)_{s\ge 0}$ denotes a standard Brownian motion in $\R^m$ (see \cite{Newman1999}). Specifically, 
	\begin{equation}\label{muL}
		\frac{d\mu}{dx} \propto \exp(-L/\eps)\,.
	\end{equation}
	In general, the invariant measure $\mu$ is not explicitely known, but it can be shown that if the system is completely controllable, then it is possible to efficiently sample it by Monte Carlo, exploiting ergodicity of the process (see \cite{Shirikyan2017,Hairer2011}). 
	
	We will come back to the question of stochastic approximations of the controllability function or some of its variants in Section \ref{sec:SOC} below. 
	
	\subsection{Finite-time controllability to a set}
	
	We consider now two modifications of (\ref{scherpenCost})--(\ref{scherpenConstraint}), the first of which is to consider a finite time horizon rather than an infinite time horizon and the second of which is to consider a target set rather than a target state. 
	
	To this end, let $B\subset\R^n$ be a closed bounded set with smooth boundary and define the cost functional
	\begin{equation}\label{setCost}
		K(u) = \frac{1}{2}\int_{0}^{T} |u(s)|^2 ds
	\end{equation}
	that we seek to minimise over controls $u$ subject to the constraint 
	having a finite value under the constraint
	\begin{equation}\label{setConstraint}
		x(t)=x\,,\quad x(T)\in B\,.
	\end{equation}
	Let $V\colon\R^n\times[0,T]\to[0,\infty]$ denote the value function 
	\begin{equation}\label{valuefct}
		V(x,t) = \min_{u\in\cA} \frac{1}{2}\int_{t}^{T} |u(s)|^2 ds + h(x(T))\,,
	\end{equation}
	where the terminal cost 
	\begin{equation}\label{setBC}
		h(x) = \begin{cases}
			0\,, & x\in B\\ \infty\,, & x\notin B
		\end{cases}\,,
	\end{equation}
	is a penalisation term that imposes the constraint $x(T)\in B$, and   
	$\cA$ denotes the set of admissible controls, such that (\ref{ode}) has a unique solution on the interval $[t,T]$. 
	It follows by the dynamic programming principle, that the value function must satisfy the Hamilton-Jacobi-Bellman (HJB) equation (see \cite{fleming2006})
	\begin{equation}\label{setHJB}
		-\frac{\partial V}{\partial t} = f\cdot\nabla_x V - \frac{1}{2}|\nabla_x V|^2_{gg^T}\,,\quad V(x,T)=h(x)\,.
	\end{equation}
	
	Conversely, by a modification of the verification theorem in \cite{Scherpen1993}, we obtain the following: 
	\begin{prop}
		If $V$ is a classical solution of the HJB equation (\ref{setHJB}), then $V$ is the value function associated with (\ref{setCost})--(\ref{setConstraint}). The value function has the form 
		\begin{equation}\label{value}
			V(x,t) = \frac{1}{2}\int_{t}^{T} |u^*(s)|^2 ds + h(x^*(T))\,,
		\end{equation}
		with 
		\begin{equation}\label{setOC}
			u^*(s) = - (g(x^*(s)))^T\nabla_x V(x^*(s),s)
		\end{equation}
		being the unique optimal control that minimises (\ref{setCost}).
	\end{prop}
	\begin{pf}
		By differentiating the solution $V\in C^{1,1}$ of (\ref{setHJB}) along $x(s)=x(s;u)$ at time $s\in(t,T)$ with $0\le t<T$, we obtain by completing the square:
		\begin{align*}
			\frac{d}{ds} V(x(s),s) & = \frac{\partial V}{\partial s} + \nabla_x V\cdot \left(f + g u\right)\\
			& = \frac{1}{2}|\nabla_x V|^2_{gg^T}  + \nabla_x V\cdot (gu) \\
			& =  - \frac{1}{2}|u|^2 + \frac{1}{2}|u + g^T\nabla_x V|^2\,.
		\end{align*} 
		Integrating the last expression from $t$ to $T$, we obtain   
		\begin{equation}
			V(x,t) - V(x(T),T) =  \frac{1}{2}\int_t^T|u|^2 -|u + g^T\nabla_x V|^2ds
		\end{equation} 
		As a consequence, 
		\begin{equation}
			V(x,t) \le \frac{1}{2}\int_t^T|u_s|^2ds + h(x(T))\,,
		\end{equation}
		where equality is attained if and only if $u = -g^T\nabla_xV$. 
	\end{pf}

	\section{A stochastic control problem}\label{sec:SOC}

	Even computing approximate solutions to the HJB equation (\ref{setHJB}) is not an easy task, especially for high-dimensional systems. Moreover, the dynamic programming equation (\ref{setHJB}) is likely not to have a classical, i.e. smooth solution, which produces additional challenges for numerical algorithms to solve optimal control problems. 
	
\begin{figure}
		\begin{center}
			\includegraphics[width=0.285\textwidth]{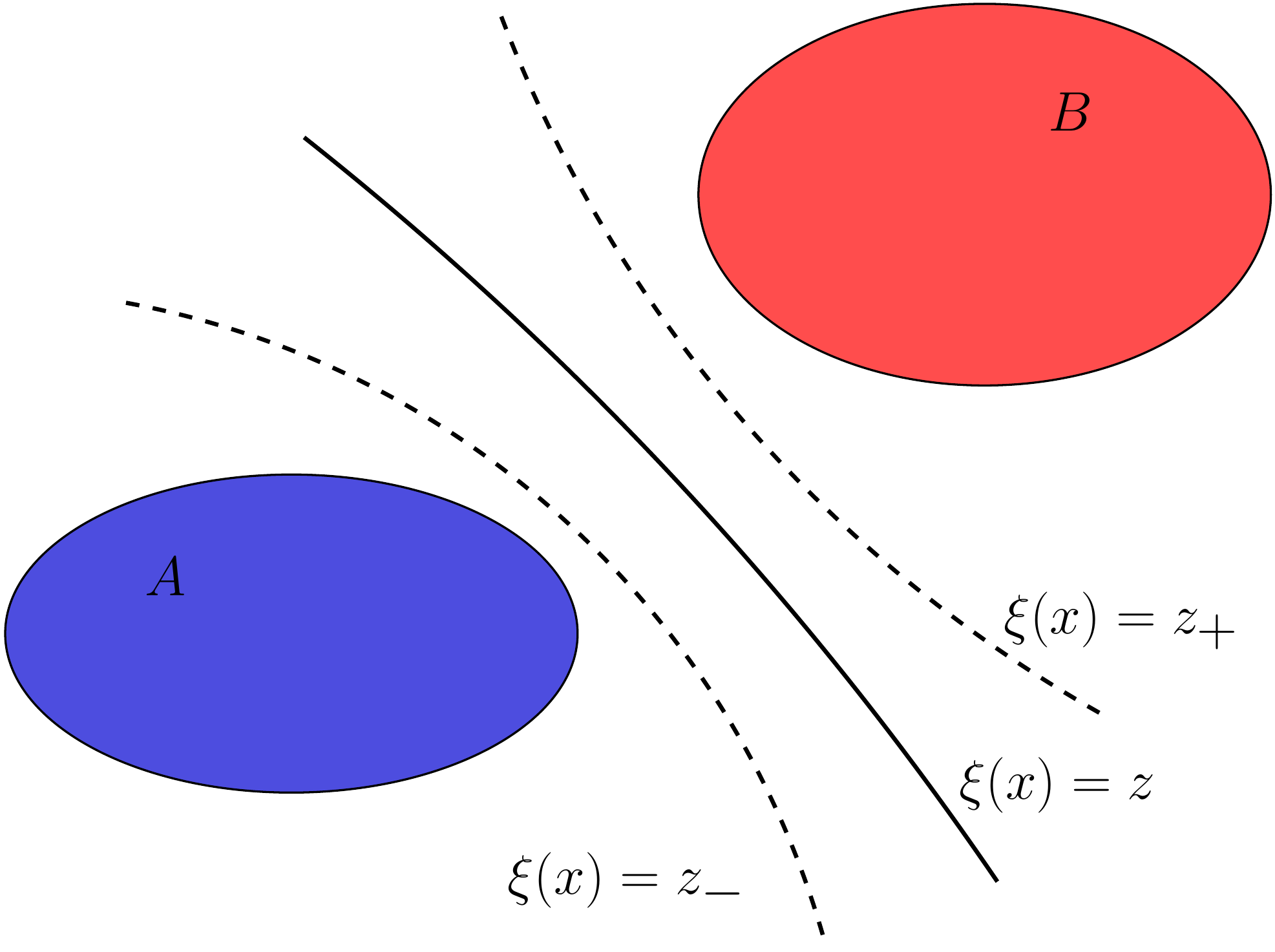}
			\caption{Reach-avoid problem with an avoidance set $A$ and a target set $B$ that are  characterised as sublevel or superlevel sets of some collective variable $\xi\,.$}\label{fig:rc}
		\end{center}
	\end{figure}

	Therefore we pick up the idea of stochastic approximations of the value function again, and suggest to add a small viscous regularisation term to (\ref{setHJB}) that mimics controller noise. To this end, we consider the controlled SDE
	\begin{equation}\label{sdeCTR}
		dX_s = f(X_s)ds + g(X_s)v_sds + \sqrt{\eps} g(X_s)dW_s\,.
	\end{equation}
	for some $\eps>0$ that is assumed to be small. Further let $A,B\subset\R^n$ denote 
	two disjoint closed bounded sets, i.e. $A\cap B=\emptyset$ as  illustrated in Figure \ref{fig:rc}. We call $C=(A \cup B)^c$ and define the random terminal time 
	\begin{equation*}
		\tau = \inf\{s>t\colon (X_{s},s)\notin C\times [0,T)\}\,, 
	\end{equation*}
	to be the stopping time that either returns $s=T$ or $s$ such that $X_{s}$ leaves the set $C$, whichever comes first.
	Now, we consider the following modification of the HJB equation (\ref{setHJB}): 
	\begin{equation}\label{setHJB2}
		-\frac{\partial V^\eps}{\partial t} = \cL^\eps V^\eps - \frac{1}{2}|\nabla_x V^\eps|^2_{gg^T}\,,\quad V^\eps|_{E^+}=h\,,
	\end{equation}
	where $E^{+} = \left(\partial C\times [0,T)\right)\cup\left(C\times \{T\}\right)$ is the terminal set of the augmented process $(X_{s},s)$, 
	and $\cL^\eps$ is a second-order differential operator given by 
	\begin{equation}\label{infGen}
		\cL^\eps = \frac{\eps}{2}gg^T\colon\nabla^2 + f\cdot\nabla\,.
	\end{equation}
	Here we use the notation $Q\colon P={\rm tr}(Q^TP)$ to denote the inner product between square matrices. 
	
	The operator $\cL^\eps$ that is defined on a suitable subspace of the space of twice continuously differentiable functions is the infinitesimal generator of the uncontrolled SDE
	\begin{equation}\label{sdeps}
		dU_s = f(U_s)ds + \sqrt{\eps} g(U_s)dW_s\,.
	\end{equation}

	Equation (\ref{setHJB2}) is the dynamic programming equation of the following optimal control problem: minimise
	\begin{equation}\label{costSDE}
		K^\eps(v) = \bE\!\left[\frac{1}{2}\int_{t}^{\tau} |v_s|^2 ds + h(X_\tau)\right]\,,
	\end{equation}
	subject to (\ref{sdeCTR}). Note that the above optimal control problem is only a minor modification of the deterministic reachability problem associated with (\ref{value}), in that we have replaced (\ref{ode}) by its weakly perturbed stochastic counterpart (\ref{sdeCTR}) and changed the definition of the terminal time from $T$ to the minimum of $T$ and the first hitting time of the boundary of $C$. 
	
	The corresponding value function $V^\eps$ can be regarded as the minimum control energy that is needed to reach a given target set $B$ before time $T$, while avoiding the set $A$. If we set $A=\emptyset$ and assume that a trajectory that hits $B$ before time $T$ remains in $B$ until time $T$, we recover essentially a noisy version of the situation in (\ref{setCost})--(\ref{setConstraint}).
	
	\begin{rem}
		The precise relation between the corresponding value functions $V^\eps$ and $V$ is addressed by the theory of viscosity solutions, for which we refer to the textbook by \cite{fleming2006} and the references therein. 
	\end{rem}
	
	\subsection{Computing the value function: hitting probabilities}
	
	Even though second-order HJB equations such as (\ref{setHJB2}) tend to have smoother solutions than their first-order counterparts, like (\ref{setHJB}), solving them numerically is notoriously difficult (see \cite{kushner2001}).
	
	It turns out, however, that the underlying nonlinear control problem belongs to a class of so-called \emph{linearly-solvable stochastic control problems} that have been studied by \cite{todorov2012}, \cite{schutte2012optimal} and others and can be solved by other means (cf.~also \cite{Nuesken2020}). The corresponding theory goes back to Fleming and co-workers, and we mention only the seminal article \cite{fleming1977}.

	We will now derive two different stochastic representations of the value function, the first of which is based on a duality argument, and the second of which uses a stochastic representation of the semilinear HJB equation (\ref{setHJB2}). 
	
	\subsubsection{Feynman-Kac representations of the value function}
	
	Assuming that $V^\eps$ is a classial solution to (\ref{setHJB2}), we can consider the logarithmic transformation 
	\begin{equation}\label{logTrafo1}
		\psi^\eps(x,t) = \exp\left(-\frac{1}{\eps}V^\eps(x,t)\right)\,.
	\end{equation}
	By chain rule,
	\begin{equation}\label{logTrafo2}
		\eps\exp\left(V^\eps/\eps\right)\cL^\eps\exp\left(-V^\eps/\eps\right) = \frac{1}{2}|\nabla V^\eps|_{gg^T} - \cL^\eps V^\eps\,,
	\end{equation}
	which, since $V^\eps$ is bounded on any compact subset that does not intersect with the terminal set $E^+$, implies that $\psi^\eps$ is the solution to the linear boundary value problem
	\begin{equation}\label{setBVP}
		\left(\frac{\partial}{\partial t} + \cL^\eps\right)\psi^\eps=0\,,\quad \psi^\eps|_{E^+}=\one_{B}\,,
	\end{equation}
	where we have used that, for all $\eps>0$, 
	\begin{equation}
		\exp(-h(x)/\eps) = \begin{cases}
			1\,, & x\in B\\ 0\,, & x\notin B
		\end{cases}\,.
	\end{equation}
	By the Feynman-Kac theorem (e.g. \cite{oksendal2003}), the solution to (\ref{setBVP}) can be recast as 
	\begin{equation}
		\psi^\eps(x,t) = \bE\!\left[\one_{B}(U_\tau)\big| U_t=x\right]\,,
	\end{equation} 
	in other words, 
	\begin{equation}
		\psi^\eps(x,t) = P\!\left(U_\tau\in B\big|U_t=x\right)\,,
	\end{equation} 
	where $U=(U_s)_{s\ge t}$ is the solution of the uncontrolled SDE (\ref{sdeps}). 
	Hence, 
	\begin{equation}\label{logTrafo3}
		P\!\left(U_\tau\in B\big| U_t=x\right) = \exp\left(-\frac{1}{\eps}V^\eps(x,t)\right)\,,
	\end{equation}
	where the resemblance with (\ref{muL}) is no coincidence. 
	
	Related reachability problems appear in reliability and have been studied in connection with obstacle avoidance problems and stochastic hybrid systems; see \cite{Bujorianu2009,mohajerin2016,Elguindy2017}.
	For high-dimensional systems, solving the linear boundary value problem (\ref{setBVP}) is by hardly any means simpler than solving the equivalent HJB equation (\ref{setHJB2}), however the Feynman-Kac representation allows for approximating the value function by Monte Carlo, especially since there is a large zoo of efficient numerical methods to compute hitting probabilities in high dimension (e.g. \cite{Elber,Wales2009}).
	

	There are two interesting limit cases associated with the stochastic representation of the value function $V^\eps$: 
	
	\paragraph*{1. Hitting probability of $B$:} If we consider the case $A=\emptyset$, then leaving the set $C$ means hitting the set $B$, and thus the value function is the negative logarithm of the probability that the uncontrolled process reaches $B$ before time $T$:  
	\begin{equation}\label{hittingProb}
		V^\eps(x,t) = -\eps \log P\!\left(\tau \le T\big| U_t=x\right) \,.
	\end{equation}
	Here $\tau=\tau_B$ is the first hitting time of $B$. 
	
	\paragraph*{2. Committor probability of hitting $B$ before $A$: } 
	Suppose that the sets $A,B$ have positive measure and can be both reached in finite time with positive probability. Now letting $T\to\infty$, it follows that $\tau$ converges with probability one to the first hitting time of $A\cup B$. Moreover, by the strong Markov property of the process $U$, the function $\psi^\eps$ becomes independent of $t$ as $T\to\infty$, so that\footnote{Since $f$ and $g$ are not explicitely time-dependent and the time horizon is infinitely large, it does not matter when to start.} 
	\begin{equation}\label{committorProb}
		V^\eps(x) = -\eps \log P\!\left(\tau_B < \tau_A \big| U_0=x\right) \,.
	\end{equation}
	Here $\tau_A$ and $\tau_B$ denote the first hitting times of the sets $A$ and $B$. The probability $ P(\tau_B < \tau_A)$ is called the \emph{committor probability} or \emph{potential} and plays a prominent role in statistical mechanics (e.g., see \cite{TPT,Bovier2015}).

	\section{Hitting probabilities for small noise systems}\label{sec:smallnoise}
	Let us consider the first hitting time of a set $B$ for the dynamics given in \eqref{sdeps} for $T\to \infty$, i.e.
	\begin{equation}
		\tau = \inf\left\{ s> t\colon  U_s \in  \bar B  \right\}\,,
	\end{equation}
	where $\bar B$ is the closure of $B$ and in the vanishing viscosity limit, i.e. $\eps \to 0$. In these limits, the first hitting time is exponentially distributed (see e.g. \cite{schutte2012optimal,neureither2017time}), so that the parameter of interest is $\bE(\tau)$, since 
	\[
	P\!\left(\tau \le T\big| U_t=x\right) \simeq 1 - \exp\left(-\frac{T-t}{\bE(\tau)}\right)\,.
	\]

	For the corresponding deterministic (or uncontrolled) dynamics $\dot u = f(u)$ we assume without loss of generality that $0 \in B^c$ is an asymptotically stable fixed point and $B^c$ is invariant. In this case, large deviations theory  provides an expression for $\bE(\tau)$ in the vanishing viscosity limit $\eps \to 0$ via the quasi-potential (see e.g. \cite{freidlin1998random}). Roughly speaking, the quasi-potential resembles the controllability function in that it measures the amount of noise needed for the process \eqref{sdeps} to reach a given state $y \in \R^n$ at time $T>0$ starting from 0. The restriction here lies in the non-degeneracy assumption for the diffusion matrix $g$, i.e. it is assumed that $gg^T >0$. The work of \cite{zabczyk1985exit} gives a natural extension of this theory by resorting to the corresponding deterministic control problem \eqref{ode}-\eqref{scherpenConstraint} so that the controllabilty function takes the role of the quasi-potential. 
	
	To be more precise, in \cite{zabczyk1985exit} it is  shown that for any $x$ close enough to $0$ it holds that 
	\begin{equation} \label{MFET:contr}
		\lim\limits_{\eps \to 0} \eps \ \log \bE\!\left(\tau \big| U_t = x\right) = \inf\limits_{y \in \partial B} L(y)\,,
	\end{equation}
	here $\partial B $ is the smooth boundary of the domain $B$ and $L$ is the solution to \eqref{scherpenHJB}--(\ref{scherpenStable}).
	
	Dealing with the case of vanishing noise, the set $B^c$ contains a subset that is metastable with respect to the stochastic dynamics (since it contains an asymptotically stable fixed point and is assumed to be invariant under the deterministic dynamics), so that Monte Carlo methods to estimate $\bE(\tau)$ will perform poorly. This makes (\ref{MFET:contr}) particularly attractive for dynamics whose controllability function $L$ is explicitely computable, and we will discuss two cases exemplarily.
	
	\subsection{Linear systems}
	
	Let us concretise the above result for the case that $f$ is a linear function, i.e. $f(x) = A x$, where $A \in \R^{n \times n}$ and $g=C \in \R^{n \times m}$ are constant matrices. The small noise assumption makes a linearisation of the dynamics reasonable. We assume that: 
	\begin{itemize}
		\item[(A)] The matrix $A$ is Hurwitz (all eigenvalues of $A$ have strictly negative real part).
		\item[(B)] The matrix pair $(A,C)$ is controllable, i.e. 
		\[\mathrm{rank}(C|AC|A^2C|\ldots|A^{n-1}C) = n \,.\]
	\end{itemize} 
	
	The following result provides an explicit and computable expression for the mean first exit time. 
	
	\begin{thm}[\cite{zabczyk1985exit}]
		Let $B=\{x \in \R^n: |x| > r \}$. Then, under the previous assumptions, 
		\[\inf\limits_{y \in \partial B} L(y)= \frac{r^2}{2\lambda_{\max}(\Sigma)}\,,\] 
		where $\Sigma>0$ is the unique solution to the Lyapunov equation 
		\[A\Sigma + \Sigma A^T = - CC^T\] 
		and $\lambda_{\max}(\Sigma)$ is the maximal eigenvalue of $\Sigma$. 
	\end{thm}

	Note that $\eps\Sigma$ is the covariance matrix (i.e. the controllability Gramian) of the associated invariant measure \[\frac{d\mu}{dx}  \propto \exp\left(-\frac{1}{2\eps}x^T \Sigma^{-1} x \right)\] which is in accordance with \eqref{muL}.   \\
	
	\subsection{Perturbed Port-Hamiltonian systems}
	Consider the stochastic counterpart of \eqref{phs}, also known as (underdamped) Langevin equation 
	\begin{equation}\label{langevin}
		dU_s = (J-D) \nabla H(U_s)\, ds + \sqrt{2\eps D} dW_s\,,
	\end{equation}
	for which the controllability function is known to be \eqref{V-H}. Here $U=(q,p) \in \R^{2d}$ is the state vector and $J=-J^T$ and $D=D^T\ge 0$ denote the constant structure and friction matrices. Typically, 
	\begin{equation}\label{specialJD}
		J =  \begin{pmatrix} 0& I_{d \times d} \\ -I_{d\times d} & 0 \end{pmatrix},\;  D =  \begin{pmatrix} 0& 0 \\ 0 & R \end{pmatrix} \,,
	\end{equation}
	with $R=R^T>0$, even though this special form of (\ref{langevin}) is not required in what follows. 
	
	As we will show, the controllability function is related to the spectral gap of the associated generator $\cL^\eps$, i.e. the first non-zero eigenvalue $\lambda_1$ of the linear operator defined in \eqref{infGen}, that is again related to (\ref{MFET:contr}). 
	
	Let us introduce some notation and make additional assumptions as the domain of interest has to be chosen carefully yet not violating the previous assumptions:
	Let 
	\begin{equation}\label{H}
		H(q,p) = \frac{1}{2} |p|^2 + \Phi(q)\,,
	\end{equation}
	where $\Phi:\R^d \to \R$ denotes the potential energy and we assume without loss of generality that $0$ is a local minimum of $\Phi$ with $\Phi(0)=0$. We further suppose that the global minimum  of $\Phi$ is attained at $q_{\min}$. Then, among all possible paths $\gamma:[0,1] \to \R^d$ connecting the two minima, such that  $\gamma(0)=0$ and $\gamma(1)=q_{\min}$, we choose the one that passes over the unique lowest energy barrier that we denote by $q^*$, i.e. 
	\[\Phi(q^*) = \inf\limits_{\gamma} \sup\limits_{t \in [0,1]} \Phi(\gamma(t))\,.\]  
	Define 
	\begin{equation} \label{phs:B}
		B = \left\{ (q,p) \in \R^{2d}: |(q,p)| > |q^*| \right\} 
	\end{equation} 
	and assume that $0 \in B^c$ is an asymptotically stable fixed point for the deterministic dynamics $\dot u = f(u)$ and that 
	\begin{equation} \label{phscondB}
		\inf\limits_{(q,p) \in \partial B} H(q,p) = \Phi(q^*)\,.
	\end{equation}
	Since $ \left\{(q,0): |q|=|q^*| \right\} \subset \partial B\,,$  \eqref{phscondB} and \eqref{H} entail that  \[ \Phi(q^*) = \inf\limits_{(q,0) \in \partial B} H(q,p) =  \inf \limits_{\left\{ q: |q|=|q^{*}| \right\}} \Phi(q)\,,\]
	so that we may equivalently require the following two conditions to hold
	\begin{itemize}
		\item[(a)] $\Phi(q) \geq \Phi(q^*) \ \forall q: |q|=|q^*|$
		\item[(b)] $H(q,p) \geq H(q^*,0) \ \forall (q,p) \in \partial B\,.$
	\end{itemize}  Note that even though the lowest energy barrier is uniquely determined by $q^*$, this does not imply that the infimum $\inf\limits_{(q,p) \in \partial B} H(q,p)$ is unique.
	\begin{prop}\label{thm:phs}
		Under the above assumptions and for small enough $\eps$ the first non-zero eigenvalue $\lambda_1$ of $\cL^\eps$ in \eqref{infGen}
		and $\bE(\tau \big| U_t=x)$ with $x$ being close enough to 0 satisfies
		\[\lambda_1 \simeq - 1/ \bE\!\left(\tau \big| U_t=x\right)\,,\]
		where $\tau$ is the first hitting time of $B$ defined in \eqref{phs:B}.
	\end{prop}
	\begin{pf}
		Following \cite{herau2008tunnel}, under the given assumptions, it holds that
		\begin{equation*}
			\lambda_1 \approx - c \exp(- \Phi(q^*)/\eps)
		\end{equation*}
		for small enough $\eps$ and  a positive constant $c = \tilde c + \mathcal{O}(\eps)\,,$ with $ \tilde c > 0\,.$ This yields 
		\begin{equation*}
			\lim\limits_{\eps \to 0} \eps \log(-\lambda_1) = - \Phi(q^*)\,.
		\end{equation*}
		On the other hand, using \eqref{MFET:contr}, it follows that 
		\begin{align*}
			\lim\limits_{\eps \to 0} \eps \log \bE\!\left(\tau \big| U_t=x\right) &= \inf\{ H(q,p) \colon |(q,p)|=|q^*|\}\\
			&= \Phi(q^*)\,.
		\end{align*}
		As a consequence, 
		\begin{equation*}
			\lim\limits_{\eps \to 0} \eps \log(-\lambda_1) = \lim\limits_{\eps \to 0} \eps \ \log (1/\bE(\tau | U_t = x))
		\end{equation*}
		which we write as $-\lambda_1 \simeq 1/\bE(\tau | U_t = x)$.
		
	\end{pf}

	\subsection{Controllability function for collective variables}
	

	Often  a target set or a target state is defined by a subset of coordinates; for example, a target set or state may involve configurations, but not momenta, or it may be defined by a collective variable, i.e. a function 
	\[\xi\colon\R^n\to \R^k\] 
	with $k<n$. Specifically, we consider the Langevin dynamics (\ref{langevin}) with $U=(q,p)\in\R^{2d}$ and assume that
	\[
	q=(\bar{q},\eta)\in\R^{l}\times\R^{n-l},\; p=(\bar{p},\zeta)\in\R^{l}\times\R^{n-l}\,.
	\]
	We call $\bar{q},\bar{p}$ the resolved variables, and $\eta,\zeta$ the unresolved variables. 
	For simplicity, we suppose that $\xi=(\xi_1,\xi_2)$, with $\xi_1(q,p)=\bar{q}$ and $\xi_2(q,p)=\bar{p}$. 
	Our aim will be to express the controllability function $L$ as a function of the resolved variables  $z=(\bar{q},\bar{p})$ only, assuming that there exists a function $\bar{L}\colon\R^{2l}\to [0,\infty)$, such that  
	\begin{equation}\label{Lbar}
		\|L - \bar{L}\circ \xi\|\approx 0
	\end{equation}
	in some appropriate sense (e.g. uniformly on any compact subset of $\R^n$ or in mean square error).

	\paragraph*{Some notation} In order to derive a coarse grained representation of the controllability function for collective variables, it is convenient to introduce the Hilbert space 
	\[
	\cH = \left\{\varphi\colon \R^{2d}\to\R\, \colon \int_{\R^{2d}} |\varphi(x)|^2\,d\mu(x)<\infty \right\}
	\]
	and endow it with the scalar product 
	\[
	\bk{\varphi}{\psi}_\mu = \int_{\R^{2d}} \varphi(x)\psi(x)\,d\mu(x)\,.
	\]
	where $\mu$ is a probability measure on $\R^{2d}$. 
	We call $\|\cdot\|_\mu$ the induced norm on $\cH$, and choose $\mu$ to be 
	\begin{equation}\label{muCan}
		\frac{d\mu}{dx} = \frac{1}{Z}\exp\left(-H(x)/\eps\right)\,,
	\end{equation}
	with $Z>0$ normalising the total probability to one. %
	Note that $\mu$ is the unique invariant measure of (\ref{langevin}).  
	
	Further let $\bE_\mu(\cdot)$ denote the expectation \wrt $\mu$ and call $\bE_\mu(\cdot|z)$ the corresponding conditional expectation for any given $z=(\bar{q},\bar{p})$. 
	By the properties of conditional expectations, the mapping 
	\begin{equation}
		\pi\colon \cH\to\bar{\cH}\subset\cH,\;\varphi \mapsto \bE_\mu(\varphi|z) =:\bar{\varphi}(z)
	\end{equation}
	is an orthogonal projection onto the (closed) subspace of $\cH$ that contains only functions of  $z$, where orthogonality is understood \wrt the $\mu$-weighted scalar product  $\bk{\varphi}{\psi}_\mu$ on $\cH$---in other words, 
	\[
	\bk{\varphi - \bar{\varphi}}{\bar{\varphi}}_\mu = 0\quad \forall \varphi\in\cH\,.
	\]
	Hence it enjoys the best approximation property:
	\[
	\|\varphi - \bar{\varphi}\|_\mu = \min_{\psi\in\bar\cH} \|\varphi - \psi\|_\mu\,.\]
	
	\begin{defn}
		The quantity 
		\begin{equation}\label{freeEnergy}
			\bar{L}(z) =  -\eps\log\int_{\R^{2d-2l}}\exp\left(-H/\eps\right)d\eta d\zeta
		\end{equation}
		is called the \emph{free energy} in the resolved variables $z=(\bar{q},\bar{p})$. 
	\end{defn}

	The following theorem characterises the coarse grained controllability function under the assumption (\ref{Lbar}).
	\begin{thm}\label{thm:cg}
		Suppose there exists a function $\bar{L}\colon\R^{2l}\to [0,\infty)$ and some sufficiently small $\delta>0$, such that  
		\[
		L - \bar{L}\circ\xi = \delta R\,,
		\]
		where $L$ is the unique classical solution of the dynamic programming equation (\ref{scherpenHJB})--(\ref{scherpenStable}) and $\|R\|_{\mu}<\infty$. Then, the free energy (\ref{freeEnergy}) 
		solves the projected HJB equation
		\begin{equation}\label{cgHJB}
			\bar{f}\cdot\nabla \bar{L} + \eps|\nabla \bar{L}|_{\bar{D}}=0\,,\quad \bar{L}(0)=0\,,
		\end{equation}
		with $\bar{f} = \pi f_z$ under the constraint that 
		\begin{equation}\label{cgStable}
			\bar{f}_{-}^* = -(\bar{f}+2\eps\bar{D}\nabla \bar{L})\,.
		\end{equation}
		has an asymptotically stable equilibrium at $0\in\R^{2l}$. Here $\bar{f} = \pi f_z$ denotes the best approximation of the $z$-component of the vector field $f=(J-D)\nabla H$ as a function of $z$, and $\bar{D}=D_{zz}$ is the corresponding friction matrix. 
	\end{thm}
	\begin{pf}
		We refrain from giving a rigorous proof and give a formal argument instead. Plugging the ansatz $L=\tilde{L}+\delta R$ for some function $\tilde{L}=\tilde{L}(z)$ into (\ref{scherpenHJB}) and equating different powers of $\delta$, we obtain to lowest order:
		\[
		f_z\cdot\nabla \tilde{L} + \eps|\nabla \tilde{L}|_{\bar{D}}=0\,,\quad \tilde{L}(0)=0\,.
		\]
		Since $\nabla\tilde{L}$ and $\bar{D}$ are independent of the unresolved variables, the equation has a nontrivial solution iff
		\[
		\pi(f_z\cdot\nabla \tilde{L}) = \bar{f}\cdot\nabla \tilde{L}= -\eps|\nabla \tilde{L}|_{\bar{D}}\,,
		\]
		where, since $J$ and $D$ are constant, 
		\begin{align*}
			\bar{f} & = (J-D)_{zz} (\pi\nabla_z H)(z)\\
			& = (\bar{J}-\bar{D}) \frac{\displaystyle \int_{R^{2d-2l}} \nabla_z H\exp(-H/\eps)\,d\eta d\zeta}{\displaystyle \int_{R^{2d-2l}}\exp(-H/\eps)\,d\eta d\zeta}\\
			& = (\bar{J}-\bar{D}) \nabla\left(-\eps\log \int_{R^{2d-2l}}\exp(-H/\eps)\,d\eta d\zeta\right)\,,
		\end{align*}
		using the shorthands $\bar{J}=J_{zz}$ and $\bar{D}=D_{zz}$. 
		The expression in the parenthesis equals $\bar{L}$, and it is easy to see that $\tilde{L}=\bar{L}$ is a  nontrivial solution of (\ref{cgHJB}) and satisfies the stability condition (\ref{cgStable}). 
	\end{pf}
	
	\begin{rem}
		Even though we did not show that $\bar{L}$ is the unique solution to the projected HJB equation, it is likely that uniqueness follows from standard arguments using the uniqueness of $L$ and the fact that the free energy inherits most of the properties of the original Hamiltonian (e.g. smoothness, boundedness from below, etc.). We refer to the textbook by \cite{fleming2006} for further details on existence and uniqueness of HJB equations. 
	\end{rem}
	
	\paragraph*{Controllability function as potential of mean force}
	Theorem \ref{thm:cg} states that, up to an additive constant, the coarse grained controllability function as a function of $\xi$ is the logarithm of the $\xi$-marginal, i.e. the density of the pushforward measure $\mu\circ\xi^{-1}$:
	\[
	\bar{L} = -\eps\log\frac{d(\mu\circ\xi^{-1})}{dz}\,,
	\]
	which, by definition, is the free energy of the system parametrised by the resolved (i.e. macroscopic) variables $z=(\bar{q},\bar{p})$. Given that the controllability function is the value function of the optimal control problem to reach the macroscopic target state $0\in\R^k$---instead of $0\in\R^n$ that involves unresolved (microscopic) variables---from another macroscopic state $z$ with minimum energy, this finding is consistent with the conventional physical interpretation of the free energy as the \emph{potential of mean force}. 
	
	We should stress, that the above argument is strongly tied to the specific form of the dynamics (\ref{langevin})  \emph{and} the coordinate choice, in that the situation that the projection of the drift generates an invariant measure that is related to the original invariant measure by the same projection is not a triviality; even if an SDE like (\ref{sde}) has a unique invariant measure with a well-defined marginal, its projection may not even have an invariant measure; see \cite{neureither2019,siads2020}.
	
	Even for a Langevin or PHS-like dynamics like (\ref{langevin}) or (\ref{phs}) a badly chosen collective variable may not lead to a proper SDE with an invariant measure given by the corresponding marginal. For example, if $\xi(q,p)=q$ is the projection onto the configuration variable for a Langevin dynamics of the form (\ref{langevin})--(\ref{H}), then $\mu\circ\xi$ will be the position marginal with density $\exp(-\Phi/\eps)$, but the projected dynamics will be of the form $\dot{q}=0$. Thus the corresponding invariant measure will be singular (i.e. a Dirac mass at $q=q(0)$).  
	
	Nevertheless the free energy $\bar{L}=\Phi$ may be informative as a quasi-potential as Proposition \ref{thm:phs} shows.

	\section{Numerical example: double pendulum}\label{sec:num}
	
	\begin{figure}
		\begin{center}
			\includegraphics[width=0.285\textwidth]{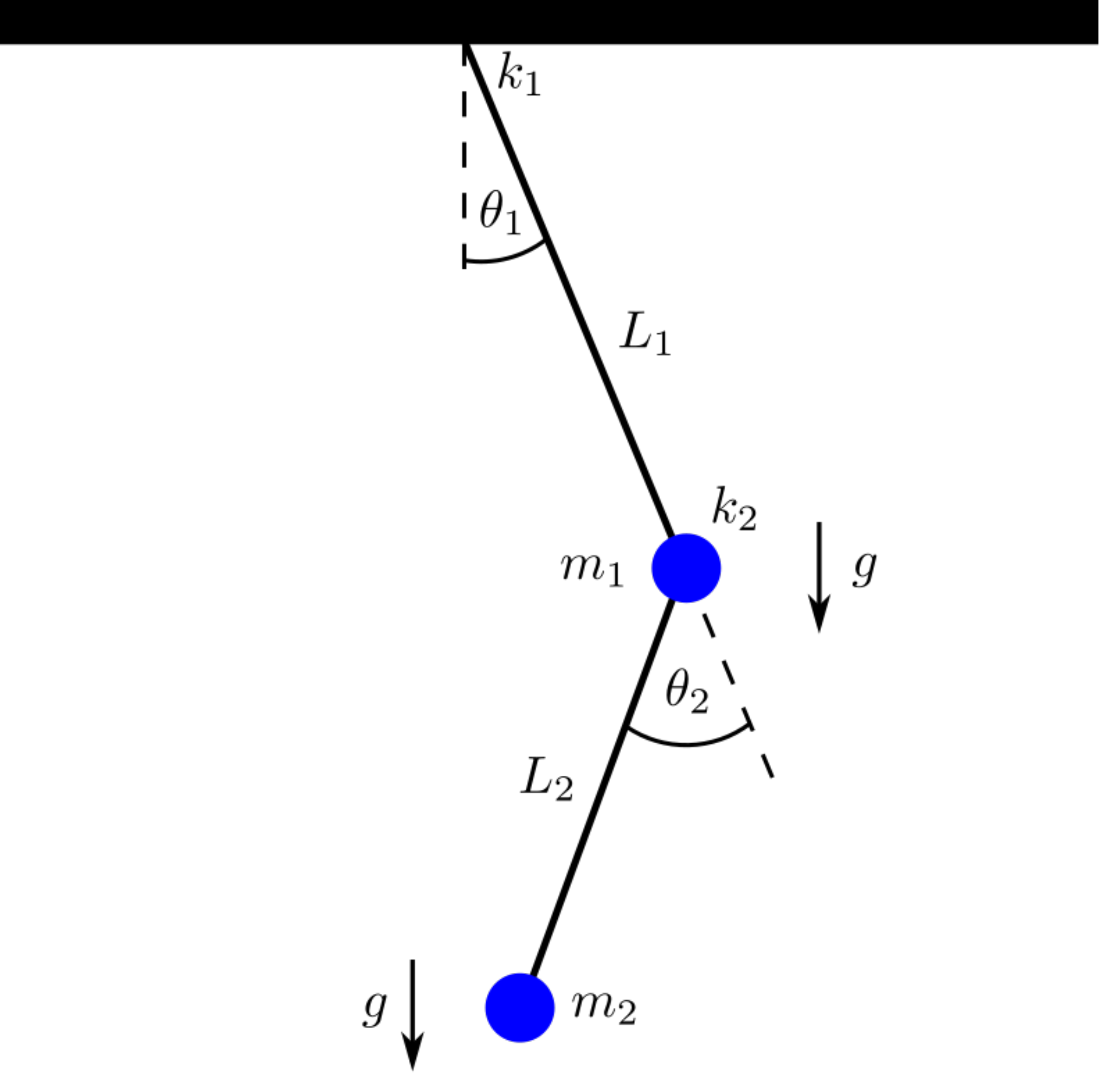}
			\caption{A double pendulum with joint angles $\theta_i$, masses $m_i$, linear torsional stiffness $k_i$, shaft lengths $L_i$, and gravitational acceleration $g$.}\label{fig:double_pendulum}
		\end{center}
	\end{figure}
	
	We consider a planar double pendulum with massless shafts (see Figure \ref{fig:double_pendulum} for details) with $m_i=k_i=L_i=1$ and $g=10$.
	In the canonical coordinates $q = (\theta_1, \theta_2)$ and $p = (p_1, p_2)$ the Hamiltonian can be written as
	\[
	H = \frac{1}{2}p^T(M(q))^{-1}p + \Phi(q)
	\]
	with the inverse mass matrix 	
	\begin{align*}
		M^{-1} = \frac{1}{2-\cos(q_2)}\begin{pmatrix}
			1 & -1-\cos(q_2)\\
			-1-\cos(q_2) & 2+3\cos(q_2)
		\end{pmatrix}
	\end{align*}
	and the potential
	\[
	\Phi = \frac{1}{2}\left(q_1^2+q_2^2\right) - 20\cos(q_1) - 10\cos(q_1+q_2)\,.
	\]
	The Langevin equation in $u=(q,p)\in\mathbb{T}^2\times\R^2$, with $\mathbb{T}^2=S^1\times S^1$ denoting the two-dimensional torus reads
	\begin{equation}\label{langevin4D}
		dU_s = (J-D) \nabla H(U_s)\, ds + \sqrt{2\eps D} dW_s\,,
	\end{equation}
	with
	\[
	J =  \begin{pmatrix} 0& I_{2 \times 2} \\ -I_{2\times 2} & 0 \end{pmatrix},\;  D =  \begin{pmatrix} 0& 0 \\ 0 & I_{2\times 2} \end{pmatrix} \,,\; 
	\]
	and the controllability function is $L = H - \Phi(0)$.
	
	\subsection{Exit from a set}	
	
\begin{figure}[b!]
		\begin{center}
			\includegraphics[width=0.5\textwidth]{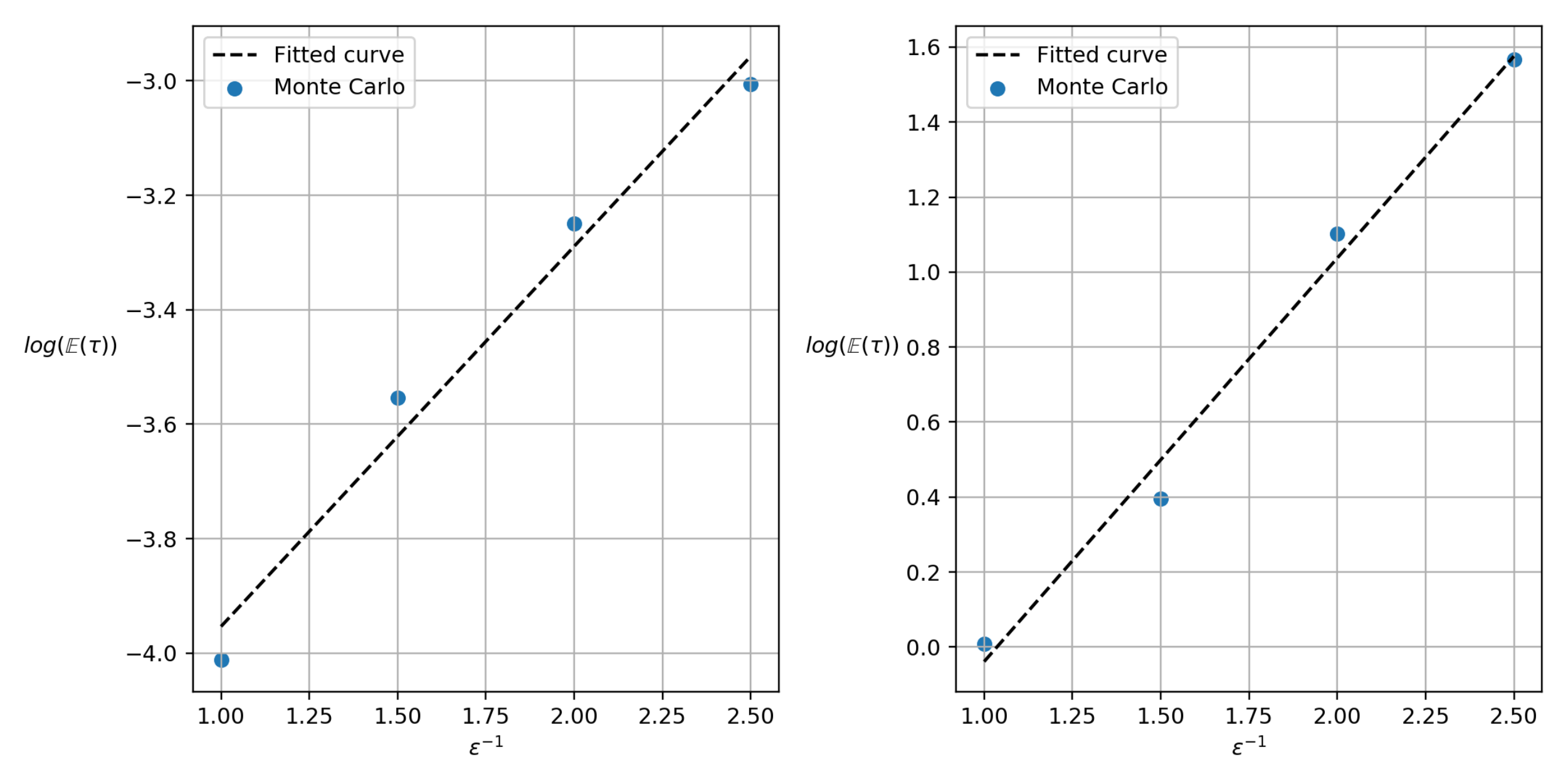}
			\caption{Logarithmic plot of the estimated mean first exit time in the non linear case as function of $1/\eps$. The left panel shows the result of the Monte Carlo simulation for the domain $\mathcal{D}_1$, the right panel for $\mathcal{D}_2$.}\label{fig:non_lin_MC}
		\end{center}
	\end{figure}

	We are interested  in reaching the boundary of one of the two sets 
	\begin{align*}
		\mathcal{D}_1 = \left\{(p,q)\in \mathbb{T}^2\times\R^2\colon|(M^{-1}p,q)| \leq 1\right\}
	\end{align*}
	or 
	\begin{align*}
		\mathcal{D}_2 = \left\{(p,q)\in \mathbb{T}^2\times\R^2\colon |q| \leq 0.5\right\}.
	\end{align*}
	The infimum is computed using the function $minimize$ from the python package scipy which yields \[\inf\limits_{(p,q) \in \partial \mathcal{D}_1} L(p,q)=0.0858 \text{ and }\inf\limits_{q \in \partial \mathcal{D}_2} L(q)=0.8539\,.\] We compare these values with the slope of the linear interpolation of $\log E(\tau)$ as a function of $1/\epsilon$ (see Figure~\ref{fig:non_lin_MC}). The mean first exit times are computed by Monte Carlo  with $120$ realisations. We obtain an estimate of $0.6642$ for the set $\mathcal{D}_1$ and $1.0761$ for $\mathcal{D}_2$.\\
	
	\subsection{Linearised systems}
	
	The linearised Langevin dynamics is obtained by replacing $H$ in (\ref{langevin4D}) by its quadratic approximation about the stable equilibrium $u_0=(0,0)$ 
	\begin{align*}
		\hat{H} = \frac{1}{2}p^T S p + \frac{1}{2}q^T N q,
	\end{align*}
	where we have omitted the additive constant $\Phi(0)$ and defined the symmetric positive definite matrices 
	\begin{align*}
		S = \begin{pmatrix}
			1 & -2\\
			-2 & 5
		\end{pmatrix}\,,\; 	N  = \begin{pmatrix}
			31 & -10\\
			-10 & 11
		\end{pmatrix}\,.
	\end{align*}
	The associated controllabillity function then is  
	\begin{align*}
		\hat{L} = \hat{H} - \Phi(0)\,.
	\end{align*} 
	
	Minimizing $\hat{L}$ on the boundary of $\mathcal{D}_2$ yields 
	\begin{align*}
		\inf\limits_{q \in \partial \mathcal{D}_2} \hat{L}(p,q) = \inf\limits_{q \in \partial \mathcal{D}_2} \frac{1}{2}q^T N q = \frac{1}{2}\lambda_{min}(N)0.5^2\,,
	\end{align*} 
	where $\lambda_{\min}(N)$ is the minimal eigenvalue of the matrix $N$.
	Observe that the Lyapunov equation corresponding to the  linearised dynamics reads
	\begin{align*}
		\begin{pmatrix}
			0 & S\\
			-N & -S
		\end{pmatrix}\Sigma + \Sigma \begin{pmatrix}
			0 & -N\\
			-S & -S
		\end{pmatrix} =  -\begin{pmatrix}
			0 & 0\\
			0 & 2I_{2\times 2}
		\end{pmatrix},
	\end{align*}
	with the unique solution 
	\[
	\Sigma = \begin{pmatrix}
		N^{-1} & 0\\
		0 & S^{-1}
	\end{pmatrix}\,.
	\]
	The minimisation problem can be solved analytically: 
	\[\inf\limits_{q \in \partial \mathcal{D}_2} \hat{L}(p,q) = 0,8572.\]
	
	\begin{rem}
		Note that -- up to additive constants -- the controllability function in the linear case associated with $\mathcal{D}_2$ is given by the free energy
		\[
		\bar{\hat{L}}(q) =  -\eps\log\int_{\R^{2}}\exp\left(-\hat{H}/\eps\right)dp\,.
		\]
		In the nonlinear case, however, this relation is only true for $\eps\to 0$:
		\begin{align*}
			\bar{L}(q) & =  -\lim_{\eps\to 0}\eps\log\int_{\R^{2}}\exp\left(-H/\eps\right)dp\\
			& = -\lim_{\eps\to 0} \left(\Phi(q) + \frac{\eps}{2}\log(\det(M(q)))\right)   = \Phi(q)\,.
		\end{align*}
	\end{rem}
	
	
	\section{Conclusions}\label{sec:sum}
	
	We have studied the set reachability of deterministic control systems and the corresponding stochastic control systems in the limit of small noise. For systems with a port-Hamiltonian structure (``Langevin dynamics''), we have derived easily computable expressions for hitting probabilities and mean first hitting times. The theoretical findings were confirmed by numerical simulations.

	\begin{ack}
		This work has been partially supported by the Collaborative Research Center \emph{Scaling Cascades in Complex Systems} (DFG-SFB 1114) through project A05, the MATH+ Cluster of Excellence (DFG-EXC 2046) through the projects EP4-4 and EF4-6, and by the StaF-EFRE project \emph{MALEDIF}. 
	\end{ack}
	
	\bibliography{control}

\end{document}